\documentclass[a4paper,12pt,oneside]{article}
\usepackage[utf8]{inputenc}
\usepackage[top=3cm, bottom=3cm, left=2.5cm, right=2.5cm]{geometry}
\usepackage{amsmath,amssymb,graphicx}
\usepackage{hyperref}
\usepackage{xcolor}
\hypersetup{
    colorlinks,
    linkcolor={blue!50!black},
    citecolor={red!40!black},
    urlcolor={green!60!black}
}

\usepackage{pgf,tikz,pgfplots}
\pgfplotsset{compat=1.14}
\usepackage{mathrsfs}
\usetikzlibrary{arrows}

\usepackage{amsmath,amssymb,graphicx,epsfig,dsfont,color,amsthm}
\usepackage{enumerate}

\def\coex{\global\advance\count1by1}\count1=0

\def\rho{\varrho}
\def\phi{\varphi}
\def\eta{\psi}
\def\qed{\nobreak\hfill $\Box$\bigskip}
\def\qedm{
\tag*{\ensuremath{\Box}}}

\newtheorem{thm}{\bfseries Theorem}
\newtheorem{lem}[thm]{\bfseries Lemma}        %% lemmas, props, cor, etc
\newtheorem{constr}{\bfseries Construction}
\newtheorem{problem}{\bfseries Problem}

\newtheorem{conj}{\bfseries Conjecture}

\newtheorem{rem}[thm]{\bfseries Remark}

\newcommand{\prf}{\noindent {\bf Proof. }}

\newcommand{\kiemel}{\textbf}

\newcommand{\red}{\ensuremath{{\mathrm{red}}}}
\newcommand{\blue}{\ensuremath{{\mathrm{blue}}}}
\newcommand{\mono}{\ensuremath{{\mathrm{mono}}}}
\newcommand{\eff}{\ensuremath{{\mathrm{eff}}}}

\definecolor{qqqqff}{rgb}{0.5,0.5,1.}
\definecolor{ffqqqq}{rgb}{0.5,0.,0.}

\definecolor{cblue}{rgb}{0.,0.,1.}
\definecolor{cred}{rgb}{1.,0.,0.}
\definecolor{cgreen}{rgb}{0.,0.4,0.2}
\definecolor{cyellow}{rgb}{1.,0.5, 0.}
\definecolor{cpurple}{rgb}{1., 0., 1.}

\author{Endre~Cs\'oka\footnote{Alfr\'ed Rényi Institute, Budapest, Hungary. The research was supported by 
ERC Synergy grant No.\ 810115.} 
\and Zolt\'an~L.~Bl\'azsik\footnote{MTA--ELTE Geometric and Algebraic Combinatorics Research Group, Budapest, Hungary. The research was supported by the Hungarian National Research, Development and Innovation Office, OTKA grant no. SNN 132625.}
\and Zolt\'an~Kir\'aly\footnote{Department of Computer Science, E\"otv\"os Lor\'and University, Budapest, Hungary. This research was partially supported by a grant (no.\ FK 132524) from the National Development Agency of Hungary, based on a source from the Research and Technology Innovation Fund, and in part by the Project no. ED\_18-1-2019-0030 (Application-specific highly reliable IT solutions) has been implemented with the support provided from the National Research, Development and Innovation Fund of Hungary, financed under the Thematic Excellence Programme funding scheme.}
\and D\'aniel~Lenger\footnote{Department of Computer Science, E\"otv\"os Lor\'and University, Budapest, Hungary.}}

\title{Upper bounds for the necklace folding problems}

\begin{document}
\maketitle

\begin{abstract}
A necklace can be considered as a cyclic list of $n$ red and $n$ blue beads in an arbitrary order, and the goal is to fold it into two and find a large 
cross-free matching of pairs of beads of different colors. 
We give a counterexample for a conjecture about the necklace folding problem, also known as the separated matching problem. The conjecture (given independently by three sets of authors) states that $\mu=\frac{2}{3}$, where $\mu$ is the ratio of the `covered' beads to the total number of beads.

We refute this conjecture by giving a construction which proves that $\mu \le 2 \nolinebreak - \nolinebreak \sqrt 2 < 0.5858$. Our construction also applies to the homogeneous model: when we are matching beads of the same color. 
Moreover, we also consider the problem where the two color classes not necessarily have the same size. 
\end{abstract}

\section{Introduction}
In the last decades, essentially the same problem known as the necklace folding or the separated matching problem appeared in many areas of mathematics. The problem has two variants which we call the heterogeneous and the homogeneous model. Consider a necklace which consists of $N=2n$ beads, $n$ red and $n$ blue ones. In both models, the aim is to find a specific folding of the necklace defined as follows. A matching $M$ consists of $|M|$ mutually disjoint pairs of beads.

In the heterogeneous model each pair consists of one red and one blue bead while in the homogeneous model each pair consists of two beads of the same color. The matched pairs will be also called matching edges.

A matching is \emph{cross-free} if no two matching edges cross each other. That is, if the two matching edges are $ab$ and $cd$, then one arc between $a$ and $b$ is disjoint from the set $\{c,d\}$ while the other arc contains this set entirely.

A \emph{secant} partitions the necklace into two continuous arcs, $A_1$ and $A_2$. A matching is \emph{secant-respecting} if for each matching edge, one end is in $A_1$ while the other end is in $A_2$. 
We call a matching \emph{proper} if it is cross-free and secant-respecting.

Let us remark here that if we drop the secant-respecting condition, then one can easily prove that there is always a cross-free matching consisting of $n$ edges in the heterogeneous (and $n-1$ edges in the homogeneous) model.

Proper matchings were called \emph{separated} matchings in   \cite{KPT1,KPT2,Viola} where the same problem was considered in a geometric setup. We have $n$ red and $n$ blue points on a circle, the matching edges are considered as segments. A matching is non-crossing if the corresponding segments are pairwise disjoint and a non-crossing matching is separated if there is a straight-line that intersects the interior of each of its segments.

Let $M$ be a proper matching. The size $|M|$ of the matching is the number of its edges. A bead is \emph{covered} if it is contained in a matched pair, the number of the covered beads is clearly $2|M|$.
Remember that a necklace consists of $N=2n$ beads, half of them is red, the other half is blue (i.e., it is balanced). For an even integer $N$, let $\mathcal{N}(N)$ denote the set of possible balanced necklaces with $N$ beads, and let $\mathcal{M}(L)$ denote the set of proper matchings for a given necklace $L$ in the heterogeneous model, and $\mu(N,L)=\max_{M\in\mathcal{M}(L)} 2|M|$, i.e., the maximum number of covered beads in a proper matching.
Moreover, let $\mu(N)=\min_{L\in\mathcal{N}(N)} \mu(N,L)$.  Thus $\mu(N)$ is the maximum number of coverable beads in the `worst' necklace.
We are interested in $\frac{\mu(N)}{N}$, the ratio of the covered beads to the total number of beads. Remark that it is the same as $|M|/n$. Finally, let $\mu=\limsup\limits_{N\rightarrow\infty} \mu(N)$. For the homogeneous model, we similarly define $\mu^{\mathrm{hom}}(N)$ and $\mu^{\mathrm{hom}}$.

\bigskip

It is trivial that there is a proper matching of size $n/2$ in any given necklace for both models. In the heterogeneous model, one can take an arbitrary secant which cuts the necklace into two arcs each of which containing $n$ beads. Since the number of blue and red beads are the same, therefore in one of the arcs there are at least as many blue beads as red ones, and in the other arc the opposite is true. Thus we can create a proper matching using the beads of the majority color from each arc. In the homogeneous model, one can take an appropriate secant for which the two arcs have the same number of blue beads. 
Then there is a proper matching of size $\lfloor n/2\rfloor$.
That is, $\mu\ge \frac12$ and $\mu^{\mathrm{hom}}\ge \frac12$.

It was very exciting that for 20 years there were no significant improvements about this lower bound, only about the additional $o(n)$ term. However, very recently Mulzer and Valtr \cite{MV} managed to improve the lower bound of 
$\mu$ to $(1/2+\varepsilon)$ for some absolute constant $\varepsilon>0$.

The story regarding the upper bound is more  diversified.
Originally only the heterogeneous model was studied.
Lyngs{\o} and Pedersen \cite{LP} in 1999 proved that 
$\mu\le 2/3$, and they conjectured that $\mu=2/3$.
Later independently Kyn\v cl, Pach and T\'oth \cite{KPT1,KPT2} and Brevier, Preissmann and Seb\H{o} \cite{Sebo} proved the same upper bound and formulated the same conjecture.
 
 \begin{conj}\cite{LP,KPT1,KPT2,Sebo}\label{main_conj}
In the heterogeneous model, there is always a proper matching of size at least $2n/3-o(n)$, i.e., $\mu=2/3$.
\end{conj}

Actually, in \cite{KPT1,KPT2} a more refined conjecture can be found. 
For a necklace $L\in\mathcal{N}(N)$, let $\mono(L)$ denote the number of maximal monochromatic arcs, i.e., there is $\mono(L)$ color changes in the necklace, or in other words, the necklace  $L$ consist of  $\mono(L)/2$ red arcs and $\mono(L)/2$ blue arcs.

 \begin{conj}[\cite{KPT1,KPT2}]\label{conj_KPT}
If we restrict ourselves to necklaces $L$ where $\mono(L)/2=k$,
then for every constant $k$ there is always a proper matching of size at least $\frac{2k-1}{3k-2}n-o(n)$ in the heterogeneous model.
\end{conj}

However, for the strict connection to Erdős problem about non-crossing alternating path, it is enough to assume that $k=o(n)$ (see below). In this case Conjecture \ref{conj_KPT} can be read as follows.

 \begin{conj}\label{main_conj_restr}
In the heterogeneous model, there is always a proper matching of size at least $2n/3-o(n)$, i.e., $\mu=2/3$, if restrict ourselves for necklaces $L\in\mathcal{N}(N)$ where $\mono(L)=o(n)$.
\end{conj}

Surprisingly, there are several connections between our problem and some interesting questions from different topics in mathematics. In the sequel, we are going to mention some of these examples as a motivation to our study. The following problem is due to Erd\H{o}s from the late 80's.

\begin{problem}
Determine or estimate the largest number $\ell = \ell(N)$ such that, for every set of $N/2$ red and $N/2$ blue points on a
circle, there exists a non-crossing alternating path consisting of $\ell$ vertices.
\end{problem}

Kyn\v cl, Pach and T\'oth \cite{KPT1,KPT2} disproved the original conjecture
of Erd\H{o}s (stating that $\ell(N)=\frac{3}{4}N+o(N)$),
and showed the following:

\begin{thm}[\cite{KPT1,KPT2}]
There exist constants $c,c'>0$ such that $\frac{1}{2}N+c\sqrt{\frac{N}{\log N}} < \ell(N) < \frac{2}{3}N + c'\sqrt{N}$.
\end{thm}

Moreover, they conjectured that the upper bound is asymptotically tight.

\begin{conj}[\cite{KPT1,KPT2}]\label{conj_alt_path}
$|\ell(N) - \frac{2}{3}N|=o(N).$
\end{conj}

Given a necklace $L\in\mathcal{N}(N)$, let $\ell(L)$ denote the maximum length of a non-crossing alternating path. They also proved the following.

\begin{thm}[\cite{KPT1,KPT2}] \label{T2}
$\ell(N)-2\cdot \mono(L)-1 \le \mu(N) \le \ell(N)$.  
\end{thm}

About the same time Abellanas et al. \cite{Abellanas} showed a very similar construction for the same upper bound. In 2010, Hajnal and M\'esz\'aros \cite{HM,Viola} improved the lower bound on $\ell(N)$ to $N/2+\Omega(\sqrt n)$, and also gave a class of configurations reaching the upper bound.

M\'esz\'aros \cite{Viola2,Viola} investigated separated matchings and found new families of constructions containing at most $\frac{4}{3}n+O(\sqrt n)$ points in any separated matching. Furthermore, she showed that if the discrepancy is at most three, then there are at least $\frac{4}{3}n$ points in the maximum separated matching.

Our main theorem (Theorem \ref{thm_best_constr}) disproves Conjecture \ref{conj_alt_path} as well by using Theorem \ref{T2}.

Interestingly, the above mentioned problems are closely related to some applied questions about the structure of proteins and some very natural questions about drawing some geometric graphs with non-crossing straight-line edges, too. In 1999, Lyngs{\o} and Pedersen \cite{LP} studied folding algorithms in the two dimensional Hydrophobic-Hydrophilic model (2D HP) for protein structure formation. They provided some approximation algorithms so that the approximation ratio depends exactly on the size of the largest proper matching in our terminology, and conjectured that there always exist a proper matching of size at least $2n/3$. 

Moreover, there are some connections between these problems and the investigation of subsequences in circular words over the binary alphabet. One can rephrase Conjecture \ref{main_conj} with  this terminology as it states that every binary circular word of length $N$ with equal number of zeros and ones has an antipalindromic linear subsequence of length at least $2N/3-o(N)$. Recently, independently from our work, M\"ullner and Ryzhikov \cite{MR,MR2}  gave a construction (which is essentially the same as our simple construction) that yields an upper bound of $2N/3+o(N)$ for both the heterogeneous and the homogeneous models (in this latter model we are looking for a palindromic linear subsequence). It seems that they were the first who studied the homogeneous model, and they made the following conjecture.

\begin{conj}[\cite{MR,MR2}]\label{conj_MR}
$\mu^{\mathrm{hom}}=2/3$.
\end{conj}

We disprove all Conjectures above. Furthermore, we improve the best known upper bound significantly by proving the following theorem.

\begin{thm}\label{thm_best_constr}
  For Construction \ref{constr_dirty}, the size of the maximum proper matching
  is at most $(2-\sqrt 2) n+o(n)$ in both the heterogeneous and the homogeneous models (i.e., $\mu \le 2-\sqrt 2$ and $\mu^{\mathrm{hom}}\le 2-\sqrt 2$).
  Moreover, Construction \ref{constr_dirty} gives an infinite series of necklaces where $\mono(L)=o(n)$.
\end{thm}

\begin{rem}
It is not obvious how this theorem disproves Conjecture \ref{conj_KPT}.
Without the details we sketch the transition.
By Theorem \ref{thm_best_constr}, there exists a specific necklace $L_1$ with $N_1$ beads where $\mu(N_1,L_1)<0.6$. Let $k=\mono(L_1)/2$. We are giving a counterexample to Conjecture \ref{conj_KPT} for this $k$, i.e., an infinite series of necklaces $L_i$ with $N_i$ beads where $\mono(L_i)=2k$ and $\mu(N_i,L_i)<0.6$. 

We get $L_i$ from $L_1$ by replacing every bead by $i$ beads of the same color. So $N_i=iN_1$, and obviously $\mono(L_i)=\mono(L_1)=2k$. Using the fact that in bipartite graphs the weight of the maximum fractional matching is the same as the weight of the maximum matching, it is not hard to prove that $\mu(N_i,L_i)=\mu(N_1,L_1)$.
\end{rem}

\begin{rem}
The problem itself, and also our construction can be defined in a measurable sense, i.e., necklace is a circle with a measurable two-coloring on its points, and for a proper matching we also require it to be measure-preserving. This is a natural generalization of the discrete problem. Although this language was very useful for finding our counterexample, we present our result in the more classical language of discrete objects. If we used the measurable definition, we may omit the terms $o(N)$ everywhere.
\end{rem}

The paper is organized as follows. In Section \ref{simple}, we are going to provide our simple construction (Construction \ref{constr_simple}) and prove that any proper matching in it has got at most $2n/3$ edges, furthermore it works for both models. We will modify the previous construction in order to further improve the upper bound and prove our main result, Theorem \ref{thm_best_constr}, in Section \ref{dirty}. Matter of fact, the improved construction (Construction \ref{constr_dirty}) also works for both models. 

Section \ref{unbalanced} is devoted to the study of the unbalanced case where one color occurs more than the other. We give new upper bounds for the case where the number of red beads is between $N/3$ and $2N/3$ by slightly changing the previous construction.
However, here we need to distinguish between the two models. Quite surprisingly, the upper bounds which we are able to acquire for the two models are the same accidentally, but the constructions are not exactly the same.

\section{The simple construction} \label{simple}

First, we give an upper bound of $2n/3$ for both $\mu$ and $\mu^{\mathrm{hom}}$, i.e., for both the heterogeneous and the homogeneous models.
After writing down this section we learned that in a recent work by 
M\"ullner and Ryzhikov \cite{MR,MR2}, the same construction has been already given. However, we decided to keep this section as a gentle introduction to our main result.

\begin{constr}\label{constr_simple}
Let $s\ge 2$ be a integer parameter, and let $n=s^{s+1}$.
The necklace consists of $s$ large arcs, each having $s^s$ blue and $s^s$ red
beads. Let $L_1,\dots,L_s$ denote the large arcs.

$L_i$ is divided into
$s^{s-i}$ red and $s^{s-i}$ blue monochromatic arcs, the colors alternates.
A monochromatic arc inside $L_i$ always consists of $s^i$ beads.
\end{constr}

For analyzing this construction, we fix an optimal proper secant-matching pair in either the homogeneous or the heterogeneous model, and denote this
optimal matching by $M$.
The secant
may split at most two large arcs, call them $L_{p}$ and $L_{r}$. If, e.g., one end of the secant is between the large arcs $L_j$ and $L_{j+1}$, then let $p=j$. We may assume that $p<r$ (if $p=r$, then every matching edge has one end in $L_p$, so $|M|\le n/s$).

Let $x$ be the first $M$-matched bead in the sequence $L_1, L_2, \dots$, and $y\in L_q$ the bead matched to $x$; if $q<p$, then set $q=p$. Thus $p\le q\le r$, and we have the following property:
any bead in $L_{p+1}\cup\dots\cup L_{q-1}$ is matched into some $L_k$ for $1\le k\le p$ and any bead in $L_{q+1}\cup\dots\cup L_{r-1}$ is matched into some $L_k$ for $r\le k\le s$.
We consider large arcs $L_p, L_q$ and $L_r$ as exceptional.

Let $M'\subseteq M$ consist of those matching edges for which no end-vertex is
in the set $L_{p}\cup L_{q}\cup L_{r}$. Our first goal is to show that $|M'|$ cannot be too large.

Let $\ell_{i,j}$ denote the $j$th monochromatic arc of $L_i$, its length is
$|\ell_{i,j}|=s^i$. If $x$ and $y$ are the first and last $M'$-matched bead in
$\ell_{i,j}$, and $M'(x)$ and $M'(y)$ are their matched partners, then we assign the
arc spanned by $M'(x)$ and $M'(y)$ to $\ell_{i,j}$, denote this arc by
$M'(\ell_{i,j})$. 
% Ehhez tényleg fontos lenne egy ábra!

We partition a subset of the beads.
Let $P(i,j)=\ell_{i,j}\cup M'(\ell_{i,j})$ for every $i$ whenever
$p<i<q$ or $r<i\le s$, and for
every $j\le 2s^{s-i}$. Clearly these $P(i,j)$ sets are pairwise disjoint. The idea is the following. For every edge of $M'$, exactly one of its end-vertex is in some $L_i$ for either $p<i<q$ or $r<i\le s$. The other end-vertex is in some $L_k$ where $k<i$. Note that the monochromatic arcs in $L_k$ are shorter than the monochromatic arcs in $L_i$, so each monochromatic arc in $M'(\ell_{i,j})$ has length at most $s^{i-1}$.

\begin{lem}
  In the set $P(i,j)$ at most $(2/3+1/s)$ fraction of the beads are
  $M'$-covered. This is true in both the heterogeneous and the homogeneous models.
\end{lem}

\prf Suppose that there are $\lambda s^i$ $M'$-covered beads in $\ell_{i,j}$.
A monochromatic arc intersecting $M'(\ell_{i,j})$ is called covered if it
contains at least one $M'$-covered vertex, and uncovered otherwise.
If we list the  monochromatic arcs intersecting $M'(\ell_{i,j})$ in the
appropriate direction, the sizes are (not strictly) monotone increasing.
The first one and the last one are covered,
and for every covered monochromatic arc (except the last one)
the next monochromatic arc is uncovered. The last monochromatic arc has length at most
$s^{i-1}$, thus we get that inside $M'(\ell_{i,j})$ the number of $M'$-covered
beads is at most the number of $M'$-uncovered beads plus $s^{i-1}$,
that is the number of $M'$-uncovered beads is at least $\lambda s^i-s^{i-1}$.
Thus in $P(i,j)$ the number of $M'$-covered beads is exactly
$2\lambda s^i$ and the number of $M'$-uncovered beads is at least
$(1-\lambda)s^i+\lambda s^i-s^{i-1}=s^i-s^{i-1}$.
So it is enough to prove that
\[
2\lambda s^i\le (2/3+1/s)((1+2\lambda)s^i-s^{i-1}),
\]
or equivalently
$6\lambda s^2\le (2s+3)((1+2\lambda)s-1)=(2+4\lambda)s^2 + (1+6\lambda)s - 3$,
i.e., $(2-2\lambda)s^2+(1+6\lambda)s\ge 3$ which is evident because $s\ge 1$ and
$0\le \lambda\le 1$.
\qed

\begin{thm}\label{thm_simple_constr}
  For Construction \ref{constr_simple}, the size of the maximum proper matching
  is at most $2n/3+o(n)$ in both the heterogeneous and the homogeneous models.
\end{thm}

\prf The subpartition defined above covers some $n'\le 2n$ beads, among them
$(2/3+1/s)n'\le(4/3+2/s)n$ are $M'$-covered. The number of beads covered
by $M\setminus M'$ can bounded by $6s^s$. So the total number of covered
beads is at most $(4/3+2/s)n+6s^s=4n/3+8s^s=4n/3+o(n)$,
as if $n$ tends to $\infty$, then $s$ also tends to $\infty$, so
$8s^s/n=8/s$ tends to zero.
\qed

\section{Proof of Theorem \ref{thm_best_constr}} \label{dirty}

We present here our main construction showing that the size of the maximum
proper matching is at most $\alpha n$ where $\alpha$ can be arbitrary close to
$2-\sqrt 2 = 0.5857\ldots<0.5858$. 

\begin{constr}\label{constr_dirty}
Let $s\ge 2$ be a integer parameter, and let $n=s^{5s+1}$.
The necklace consists of $s$ large arcs, each having $s^{5s}$ blue and $s^{5s}$ red
beads. Let $L_1,\dots,L_s$ denote the large arcs.

$L_i$ is divided into $s^{2s-i}$ red and $s^{2s-i}$ blue arcs, the colors
alternates. Let $\ell_{i,j}$ denote the $j$\textsuperscript{th} arc of $L_i$,
where $1\le j\le 2s^{2s-i}$. The arc $\ell_{i,j}$ always consists of
$s^{3s+i}$ beads. In the next step, we will change the color of some beads in
each $\ell_{i,j}$ in the following way. Let $\lambda\le\frac{1}{2}$ be a
positive parameter and for a fixed $i$, let's divide each $\ell_{i,j}$ into
$s^{s+2i}$ intervals of size $s^{2s-i}$ and in each of these tiny intervals,
change the color of $\left \lfloor \lambda s^{2s-i} \right \rfloor$ beads
backwards from the clockwise end of the tiny interval. We will refer to those
beads whose color were changed as dust in $\ell_{i,j}$.  \emph{(See Figure \ref{lij}).}

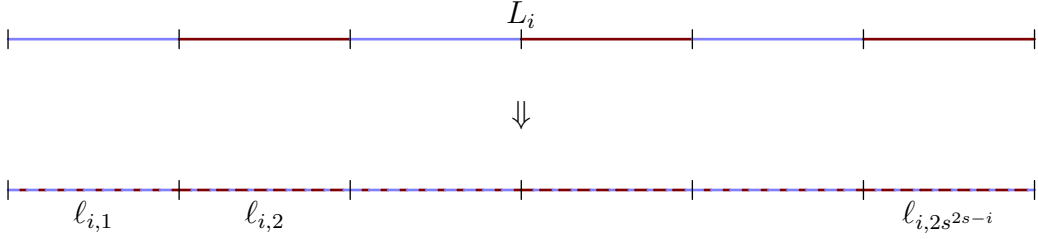
\begin{figure}[!ht]
\centering
\begin{tikzpicture}[line cap=round,line join=round,>=triangle 45,x=1.0cm,y=1.0cm, scale=0.25]

\foreach \i in {0,...,2}
{\draw [line width=1.pt,color=qqqqff] (0.+\i*18,0.)-- (9.+\i*18,0.);
\draw [line width=1.pt,color=ffqqqq] (9.+\i*18,0.)-- (18.+\i*18,0.);}

\draw (27,0) node[anchor=south] {$L_i$};
\draw (27,-4) node {$\Downarrow$};
\draw (4.5,-8) node[anchor=north] {$\ell_{i,1}$};
\draw (13.5,-8) node[anchor=north] {$\ell_{i,2}$};
\draw (49.5,-8) node[anchor=north] {$\ell_{i,2s^{2s-i}}$};

\newcommand{\arany}{0.66}
\foreach \i in {0,...,2}
{
    \foreach \j in {0,...,8}
    {
        \draw [line width=1.pt,color=qqqqff] (0.+\j+\i*18,-8.) -- (\arany+\j+\i*18,-8.);
        \draw [line width=1.pt,color=ffqqqq] (\arany+\j+\i*18,-8.) -- (1.+\j+\i*18,-8.);
        \draw [line width=1.pt,color=ffqqqq] (9.+\j+\i*18,-8.) -- (9.+\arany+\j+\i*18,-8.);
        \draw [line width=1.pt,color=qqqqff] (9.+\arany+\j+\i*18,-8.) -- (10.+\j+\i*18,-8.);
    }
    
    %\draw [line width=1.pt,color=qqqqff] (0.+\i*18,-8.)-- (9.+\i*18,-8.);
    %\draw [line width=1.pt,color=ffqqqq] (9.+\i*18,-8.)-- (18.+\i*18,-8.);
}

\foreach \i in {0,...,6}
{\draw [line width=0.5pt] (0.+\i*9,-0.5)-- (0.+\i*9,0.5);
\draw [line width=0.5pt] (0.+\i*9,-8.5)-- (0.+\i*9,-7.5);}

\end{tikzpicture}
\caption{The intervals $L_i$ and $\ell_{i,j}$}\label{lij}
\end{figure}
\end{constr}

First we bound $\mono(L)$ for a necklace $L$ given by this construction in order to prove the last statement of Theorem \ref{thm_best_constr}. 
\[
\mono(L)=\sum_{i=1}^s \sum_{j=1}^{2s^{2s-i}} 2s^{s+2i}=
4s^{3s}\sum_{i=1}^s s^i\le 4\frac{s^{s+1}-1}{s-1}s^{3s}\le 8s^{4s},
\]
which is $O(n^{4/5})=o(N)$.

We will see that for $\lambda=1 - \frac{1}{\sqrt 2}$, as $s$ tends to $\infty$, we will get the desired bound, i.e., the upper bound on the size of the proper matching tends to $2-\sqrt 2$. We will use the little-o notation, e.g., $n/s=o(n)$.

For analyzing this construction we fix an optimal proper secant-matching pair in either the homogeneous or the heterogeneous model, and denote this
optimal matching by $M$.
The secant
may split at most two large arcs, call them $L_{p}$ and $L_{r}$. If, e.g., one end of the secant is between the large arcs $L_j$ and $L_{j+1}$, then let $p=j$. We may assume that $p<r$ (if $p=r$, then every matching edge has one end in $L_p$, so $|M|\le n/s=o(n)$).

Let $M'\subseteq M$ consist of those matching edges for which no end-vertex is
inside the set $L_{p}\cup L_{r}$. Obviously, $|M|\le |M'|+2n/s=|M'|+o(n)$. We call a pair of indices $(g,h)$ \emph{bonded}, if there exists at least one edge of $M'$ connecting $L_g$ and $L_h$.
\begin{lem}
The number of the bonded pairs is at most $s-3$.  
\end{lem}

\prf
Consider the auxiliary graph with vertices $\{1, 2, \dots s\}\setminus\{p, r\}$, and connect two vertices if the corresponding pair is bonded. There is no cycle in this graph, otherwise a cycle yields a crossing in $M$. Thus the auxiliary graph may have at most $(s-2)-1$ edges.
\qed

Let $I$ be an interval. If $x$ and $y$ are the first and last $M'$-matched bead in
$I$, and $M'(x)$ and $M'(y)$ are their matched partners, then we assign the
arc spanned by $M'(x)$ and $M'(y)$ to $I$, denote this arc by $M'(I)$.

Let $(g,h)$ be a bonded pair, where $g<h$. For $i\in\{1,\dots, s^{2s-g}\}$ let's call a pair of intervals $(\ell_{g,2i-1}, \ell_{g,2i})$ \emph{$(g,h)$-regular}, if there exists $j\in\{1, \dots 2s^{2s-h}\}$ such that  $M'(\ell_{g,2i-1}\cup\ell_{g,2i})\subseteq\ell_{h,j}$. Let's call an edge of $M'$ \emph{regular} if one of the end-vertices are in a $(g,h)$-regular pair for some $g<h$. Denote the set of regular edges by $M''$. An edge of $M'$ is called a $(g,h)$-edge if its end-vertices are in $L_g$ and in $L_h$, respectively; and irregular if it is not regular. 

\begin{lem}\label{reg}
$|M'|\le |M''|+6n/s= |M''|+o(n)$.
\end{lem}

\prf
Consider a bonded pair $(g,h)$ for some $g<h$.
First we are going to bound the number of irregular $(g,h)$-edges.

Take an $i\in\{1,\dots, s^{2s-g}\}$ for which $M'(\ell_{g,2i-1}\cup\ell_{g,2i})\cap L_{h}\neq\emptyset$ but $(\ell_{g,2i-1}, \ell_{g,2i})$ is not $(g,h)$-regular. It means that either there exist a `bad index' $j\in \{1,\dots,2s^{2s-h}-1\}$ such that both $M'(\ell_{g,2i-1}\cup\ell_{g,2i})\cap\ell_{h,j}$ and $M'(\ell_{g,2i-1}\cup\ell_{g,2i})\cap\ell_{h,j+1}$ are non-empty.
We also call $j=0$ bad if $M'(\ell_{g,2i-1}\cup\ell_{g,2i})\cap L_{h-1}$ is non-empty (where $L_0=L_{s}$), and $j=2s^{2s-h}$ bad if $M'(\ell_{g,2i-1}\cup\ell_{g,2i})\cap L_{h+1}$ is non-empty
(where $L_{s+1}=L_1$).

Moreover, any $j$ can be a bad index at most once therefore the number of such $i$'s is at most $2s^{2s-h}+1$. Even if all the beads in these non-$(g,h)$-regular pairs are $M'$-covered, we got rid of at most $(2s^{2s-h}+1)\cdot 2|\ell_{g,2i}|=(2s^{2s-h}+1)\cdot (2s^{3s+g}) < 6s^{5s-(h-g)}\le 6s^{5s-1}$, since $h>g$. This is true for any bonded pairs, hence altogether we lost at most $(s-3)6s^{5s-1} < 6s^{5s}=6n/s=o(n)$ $M'$-edges.
\qed

From now on, we estimate the number of regular edges. For the sake of simplicity, we will omit the floor and ceiling functions, because the difference in the result is again $o(n)$.

For $g<h$, let's fix a bonded pair $(g,h)$. Consider a $(g,h)$-regular pair of intervals $(\ell_{g,2i-1}, \ell_{g,2i})$ such that $M'(\ell_{g,2i-1}\cup\ell_{g,2i})\subseteq \ell_{h,j}$.

Until this point, there was no difference between the homogeneous and the heterogeneous case. In the sequel, there still won't be any significant difference, the calculations work the same way in both cases. We now present the calculation for the homogeneous case, and we will assume that the `main' color of $\ell_{h,j}$ is blue (i.e. the dust is red), the `main' color of $\ell_{g,2i-1}$ is red, thus the `main' color of $\ell_{g,2i}$ is blue.

Let's denote the \emph{efficiency} of the matching $M''$ on an interval $I$ with 
$$\eff(I)=\frac{\text{\# of~} M''-\text{covered beads in } I\cup M'(I)}{|I\cup M'(I)|}.$$ 
In the following lemma, we will show that the efficiency cannot exceed $2-\sqrt 2+o(1)$ for a suitable $\lambda$.

\begin{lem}\label{eff}
$\eff(\ell_{g,2i-1}\cup\ell_{g,2i})\le 2-\sqrt{2}+o(1)$ if $\lambda=1 - \frac{1}{\sqrt{2}}$.
\end{lem}
\prf
Recall that $\ell_{g,2i-1}$ (and also $\ell_{g,2i}$) is divided into $s^{s+2g}$ red and $s^{s+2g}$ blue monochromatic intervals. We will call them $\ell_{g,2i-1}^{(\red,1)},\dots \ell_{g,2i-1}^{(\red,s^{s+2g})}, \ell_{g,2i-1}^{(\blue,1)},\dots \ell_{g,2i-1}^{(\blue,s^{s+2g})}$, where the $red$ and $blue$ indicates the color of the interval. 

Also recall that, $\ell_{h,j}$ is divided into blue and red monochromatic intervals (the color alternates) whose sizes are $(1-\lambda) s^{2s-h}$ and $\lambda s^{2s-h}$, respectively. We define numbers $a_k, b_k, c_k$ and $d_k$ for $1\le k\le s^{s+2g}$ in the following way.

Assume that the number of $M''$-covered beads from $\ell_{g,2i-1}^{(\red,k)}$ is $x$. Then \hbox{$a_k=\frac{x}{\lambda s^{2s-h}}$}, i.e., the necessary number of small red intervals (dust) from $\ell_{h,j}$ to cover that many beads.
    
Similarly, assume that the number of $M''$-covered beads from $\ell_{g,2i-1}^{(\blue,k)}$ is $x$. Then \hbox{$b_k=\frac{x}{(1-\lambda) s^{2s-h}}$}, i.e., the necessary number of blue intervals from $\ell_{h,j}$ to cover that many beads.

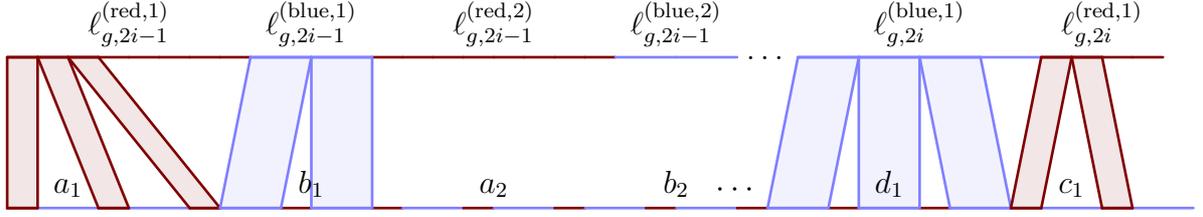
\begin{figure}[!ht]
\centering
\begin{tikzpicture}[line cap=round,line join=round,>=triangle 45,x=1.0cm,y=1.0cm, scale=0.4]
%\clip(-4.691541715971205,-11.31276094743649) rectangle (43.03902813856292,20.674592463323517);
\fill[line width=1.pt,color=ffqqqq,fill=ffqqqq,fill opacity=0.10000000149011612] (0.,5.) -- (0.,0.) -- (1.,0.) -- (1.,5.) -- cycle;
\fill[line width=1.pt,color=ffqqqq,fill=ffqqqq,fill opacity=0.10000000149011612] (1.,5.) -- (3.,0.) -- (4.,0.) -- (2.,5.) -- cycle;
\fill[line width=1.pt,color=ffqqqq,fill=ffqqqq,fill opacity=0.10000000149011612] (2.,5.) -- (6.,0.) -- (7.,0.) -- (3.,5.) -- cycle;
\fill[line width=1.pt,color=qqqqff,fill=qqqqff,fill opacity=0.10000000149011612] (8.,5.) -- (7.,0.) -- (9.,0.) -- (10.,5.) -- cycle;
\fill[line width=1.pt,color=qqqqff,fill=qqqqff,fill opacity=0.10000000149011612] (10.,5.) -- (10.,0.) -- (12.,0.) -- (12.,5.) -- cycle;
\fill[line width=1.pt,color=qqqqff,fill=qqqqff,fill opacity=0.10000000149011612] (26.,5.) -- (25.,0.) -- (27.,0.) -- (28.,5.) -- cycle;
\fill[line width=1.pt,color=qqqqff,fill=qqqqff,fill opacity=0.10000000149011612] (28.,5.) -- (28.,0.) -- (30.,0.) -- (30.,5.) -- cycle;
\fill[line width=1.pt,color=qqqqff,fill=qqqqff,fill opacity=0.10000000149011612] (30.,5.) -- (31.,0.) -- (33.,0.) -- (32.,5.) -- cycle;
\fill[line width=1.pt,color=ffqqqq,fill=ffqqqq,fill opacity=0.10000000149011612] (34.,5.) -- (33.,0.) -- (34.,0.) -- (35.,5.) -- cycle;
\fill[line width=1.pt,color=ffqqqq,fill=ffqqqq,fill opacity=0.10000000149011612] (35.,5.) -- (36.,0.) -- (37.,0.) -- (36.,5.) -- cycle;
\draw [line width=1.pt,color=ffqqqq] (0.,0.)-- (1.,0.);
\draw [line width=1.pt,color=qqqqff] (1.,0.)-- (3.,0.);
\draw [line width=1.pt,color=ffqqqq] (3.,0.)-- (4.,0.);
\draw [line width=1.pt,color=qqqqff] (4.,0.)-- (6.,0.);
\draw [line width=1.pt,color=ffqqqq] (6.,0.)-- (7.,0.);
\draw [line width=1.pt,color=qqqqff] (7.,0.)-- (9.,0.);
\draw [line width=1.pt,color=ffqqqq] (9.,0.)-- (10.,0.);
\draw [line width=1.pt,color=qqqqff] (10.,0.)-- (12.,0.);
\draw [line width=1.pt,color=ffqqqq] (12.,0.)-- (13.,0.);
\draw [line width=1.pt,color=qqqqff] (13.,0.)-- (15.,0.);
\draw [line width=1.pt,color=ffqqqq] (15.,0.)-- (16.,0.);
\draw [line width=1.pt,color=qqqqff] (16.,0.)-- (18.,0.);
\draw [line width=1.pt,color=ffqqqq] (18.,0.)-- (19.,0.);
\draw [line width=1.pt,color=qqqqff] (19.,0.)-- (21.,0.);
\draw [line width=1.pt,color=ffqqqq] (0.,5.)-- (1.,5.);
\draw [line width=1.pt,color=ffqqqq] (1.,5.)-- (2.,5.);
\draw [line width=1.pt,color=ffqqqq] (2.,5.)-- (3.,5.);
\draw [line width=1.pt,color=ffqqqq] (3.,5.)-- (4.,5.);
\draw [line width=1.pt,color=ffqqqq] (21.,0.)-- (22.,0.);
\draw [line width=1.pt,color=qqqqff] (22.,0.)-- (24.,0.);
\draw [line width=1.pt,color=ffqqqq] (24.,0.)-- (25.,0.);
\draw [line width=1.pt,color=qqqqff] (25.,0.)-- (27.,0.);
\draw [line width=1.pt,color=ffqqqq] (4.,5.)-- (5.,5.);
\draw [line width=1.pt,color=ffqqqq] (5.,5.)-- (6.,5.);
\draw [line width=1.pt,color=ffqqqq] (6.,5.)-- (7.,5.);
\draw [line width=1.pt,color=ffqqqq] (7.,5.)-- (8.,5.);
\draw [line width=1.pt,color=qqqqff] (8.,5.)-- (10.,5.);
\draw [line width=1.pt,color=qqqqff] (10.,5.)-- (12.,5.);
\draw [line width=1.pt,color=ffqqqq] (12.,5.)-- (13.,5.);
\draw [line width=1.pt,color=ffqqqq] (13.,5.)-- (14.,5.);
\draw [line width=1.pt,color=ffqqqq] (14.,5.)-- (15.,5.);
\draw [line width=1.pt,color=ffqqqq] (15.,5.)-- (16.,5.);
\draw [line width=1.pt,color=ffqqqq] (16.,5.)-- (17.,5.);
\draw [line width=1.pt,color=ffqqqq] (17.,5.)-- (18.,5.);
\draw [line width=1.pt,color=ffqqqq] (18.,5.)-- (19.,5.);
\draw [line width=1.pt,color=ffqqqq] (19.,5.)-- (20.,5.);
\draw [line width=1.pt,color=qqqqff] (20.,5.)-- (22.,5.);
\draw [line width=1.pt,color=qqqqff] (22.,5.)-- (24.,5.);
\draw [line width=1.pt,color=ffqqqq] (27.,0.)-- (28.,0.);
\draw [line width=1.pt,color=qqqqff] (28.,0.)-- (30.,0.);
\draw [line width=1.pt,color=ffqqqq] (30.,0.)-- (31.,0.);
\draw [line width=1.pt,color=qqqqff] (31.,0.)-- (33.,0.);
\draw [line width=1.pt,color=ffqqqq] (33.,0.)-- (34.,0.);
\draw [line width=1.pt,color=qqqqff] (34.,0.)-- (36.,0.);
\draw [line width=1.pt,color=ffqqqq] (36.,0.)-- (37.,0.);
\draw [line width=1.pt,color=qqqqff] (37.,0.)-- (39.,0.);
\draw (25.,5.) node[anchor=center] {$\dots$};
\draw [line width=1.pt,color=ffqqqq] (0.,5.)-- (0.,0.);
\draw [line width=1.pt,color=ffqqqq] (0.,0.)-- (1.,0.);
\draw [line width=1.pt,color=ffqqqq] (1.,0.)-- (1.,5.);
\draw [line width=1.pt,color=ffqqqq] (1.,5.)-- (0.,5.);
\draw [line width=1.pt,color=ffqqqq] (1.,5.)-- (3.,0.);
\draw [line width=1.pt,color=ffqqqq] (3.,0.)-- (4.,0.);
\draw [line width=1.pt,color=ffqqqq] (4.,0.)-- (2.,5.);
\draw [line width=1.pt,color=ffqqqq] (2.,5.)-- (1.,5.);
\draw [line width=1.pt,color=ffqqqq] (2.,5.)-- (6.,0.);
\draw [line width=1.pt,color=ffqqqq] (6.,0.)-- (7.,0.);
\draw [line width=1.pt,color=ffqqqq] (7.,0.)-- (3.,5.);
\draw [line width=1.pt,color=ffqqqq] (3.,5.)-- (2.,5.);
\draw [line width=1.pt,color=qqqqff] (8.,5.)-- (7.,0.);
\draw [line width=1.pt,color=qqqqff] (7.,0.)-- (9.,0.);
\draw [line width=1.pt,color=qqqqff] (9.,0.)-- (10.,5.);
\draw [line width=1.pt,color=qqqqff] (10.,5.)-- (8.,5.);
\draw [line width=1.pt,color=qqqqff] (10.,5.)-- (10.,0.);
\draw [line width=1.pt,color=qqqqff] (10.,0.)-- (12.,0.);
\draw [line width=1.pt,color=qqqqff] (12.,0.)-- (12.,5.);
\draw [line width=1.pt,color=qqqqff] (12.,5.)-- (10.,5.);
\draw [line width=1.pt,color=qqqqff] (26.,5.)-- (28.,5.);
\draw [line width=1.pt,color=qqqqff] (28.,5.)-- (30.,5.);
\draw [line width=1.pt,color=qqqqff] (30.,5.)-- (32.,5.);
\draw [line width=1.pt,color=qqqqff] (32.,5.)-- (34.,5.);
\draw [line width=1.pt,color=ffqqqq] (34.,5.)-- (35.,5.);
\draw [line width=1.pt,color=ffqqqq] (35.,5.)-- (36.,5.);
\draw [line width=1.pt,color=ffqqqq] (36.,5.)-- (37.,5.);
\draw [line width=1.pt,color=ffqqqq] (37.,5.)-- (38.,5.);
\draw [line width=1.pt,color=qqqqff] (26.,5.)-- (25.,0.);
\draw [line width=1.pt,color=qqqqff] (25.,0.)-- (27.,0.);
\draw [line width=1.pt,color=qqqqff] (27.,0.)-- (28.,5.);
\draw [line width=1.pt,color=qqqqff] (28.,5.)-- (26.,5.);
\draw [line width=1.pt,color=qqqqff] (28.,5.)-- (28.,0.);
\draw [line width=1.pt,color=qqqqff] (28.,0.)-- (30.,0.);
\draw [line width=1.pt,color=qqqqff] (30.,0.)-- (30.,5.);
\draw [line width=1.pt,color=qqqqff] (30.,5.)-- (28.,5.);
\draw [line width=1.pt,color=qqqqff] (30.,5.)-- (31.,0.);
\draw [line width=1.pt,color=qqqqff] (31.,0.)-- (33.,0.);
\draw [line width=1.pt,color=qqqqff] (33.,0.)-- (32.,5.);
\draw [line width=1.pt,color=qqqqff] (32.,5.)-- (30.,5.);
\draw [line width=1.pt,color=ffqqqq] (34.,5.)-- (33.,0.);
\draw [line width=1.pt,color=ffqqqq] (33.,0.)-- (34.,0.);
\draw [line width=1.pt,color=ffqqqq] (34.,0.)-- (35.,5.);
\draw [line width=1.pt,color=ffqqqq] (35.,5.)-- (34.,5.);
\draw [line width=1.pt,color=ffqqqq] (35.,5.)-- (36.,0.);
\draw [line width=1.pt,color=ffqqqq] (36.,0.)-- (37.,0.);
\draw [line width=1.pt,color=ffqqqq] (37.,0.)-- (36.,5.);
\draw [line width=1.pt,color=ffqqqq] (36.,5.)-- (35.,5.);
\draw (2.,0.) node[anchor=south] {$a_1$};
\draw (10.,0.) node[anchor=south] {$b_1$};
\draw (16.,0.) node[anchor=south] {$a_2$};
\draw (22.,0.) node[anchor=south] {$b_2$};
\draw (24.,0.) node[anchor=south] {$\cdots$};
\draw (29.,0.) node[anchor=south] {$d_1$};
\draw (35.,0.) node[anchor=south] {$c_1$};

\draw (4.,5.) node[anchor=south] {$\ell_{g,2i-1}^{(\red,1)}$};
\draw (10.,5.) node[anchor=south] {$\ell_{g,2i-1}^{(\blue,1)}$};
\draw (16.,5.) node[anchor=south] {$\ell_{g,2i-1}^{(\red,2)}$};
\draw (22.,5.) node[anchor=south] {$\ell_{g,2i-1}^{(\blue,2)}$};;
\draw (30.,5.) node[anchor=south] {$\ell_{g,2i}^{(\blue,1)}$};;
\draw (36.,5.) node[anchor=south] {$\ell_{g,2i}^{(\red,1)}$};;
\begin{scriptsize}
\end{scriptsize}
\end{tikzpicture}
\caption{The definition of $a_k, b_k, c_k$ and $d_k$}\label{fig2}
\end{figure}

We define $c_k$ and $d_k$ in the same way for $\ell_{g,2i}^{(\red,k)}$ and $\ell_{g,2i}^{(\blue,k)}$, respectively. (See Figure \ref{fig2}.) It is easy to see that $a_k\le \frac{1-\lambda}{\lambda}s^{h-g}$, $b_k \le \frac{\lambda}{1-\lambda}s^{h-g}$, and $c_k, d_k \le s^{h-g}$.

\smallskip

The number of $M''$-covered beads in $(\ell_{g,2i-1}\cup\ell_{g,2i})\cup M'(\ell_{g,2i-1}\cup\ell_{g,2i})$ is
\[
2\left(\sum a_k \lambda s^{2s-h}+\sum b_k (1-\lambda) s^{2s-h}+\sum c_k
\lambda s^{2s-h}+\sum d_k (1-\lambda) s^{2s-h}\right).
\]

In all of the sums, $k$ runs from 1 up to $s^{s+2g}$, so we use the following shorthands. Let 
$$
A=\sum_{k=1}^{s^{s+2g}}a_k, \;
B=\sum_{k=1}^{s^{s+2g}}b_k, \;
C=\sum_{k=1}^{s^{s+2g}}c_k, \;
D=\sum_{k=1}^{s^{s+2g}}d_k. 
$$

Obviously, $|\ell_{g,2i-1}\cup\ell_{g,2i}|=2s^{3s+g}$. We will give a lower bound on $|M'(\ell_{g,2i-1}\cup\ell_{g,2i})|$. To cover the $M''$-matched beads in $\ell_{g,2i-1}^{(\red,k)}$, we need at least $a_k$ monochromatic red interval from $\ell_{h,j}$, thus at least $a_k-1$ monochromatic blue interval remained unused, so $M'(\ell_{g,2i-1}^{(\red,k)})\ge (a_k-1)s^{2s-h}$. Similarly $M'(\ell_{g,2i-1}^{(\blue,k)})\ge (b_k-1)s^{2s-h}$, $M'(\ell_{g,2i}^{(\red,k)})\ge (c_k-1)s^{2s-h}$ and $M'(\ell_{g,2i}^{(\blue,k)})\ge (d_k-1)s^{2s-h}$. 

Altogether, we get that $\eff(\ell_{g,2i-1}\cup\ell_{g,2i})$ is at most 
$$2\frac{\lambda A s^{2s-h}+(1-\lambda) B s^{2s-h}+\lambda C s^{2s-h}+(1-\lambda) D  s^{2s-h}} {\left[\sum (a_k\!-\!1) s^{2s-h}+\sum (b_k\!-\!1) s^{2s-h}+\sum (c_k\!-\!1) s^{2s-h}+\sum (d_k\!-\!1) s^{2s-h}\right]+2s^{3s+g}}=$$
$$=2\frac{\lambda A+(1-\lambda) B+\lambda C+(1-\lambda) D}{[A+B+C+D-4s^{s+2g}]+2s^{3s+g}/s^{2s-h}}=$$
$$2\frac{\lambda A+(1-\lambda) B+\lambda C+(1-\lambda) D}{A+B+C+D-4s^{s+2g}+2s^{s+g+h}}.$$

Let's denote this last expression by $\eff\left(A, B, C, D\right).$ First, we will show that this expression is monotone increasing in $B$ and $D$.

$$\eff\left(A, B, C, D\right)=2(1-\lambda)+2\frac{(2\lambda-1)(A +C) - (1-\lambda)(2s^{s+g+h}-4s^{s+2g})}{A+B+C +D+(2s^{s+g+h}-4s^{s+2g})}.$$

As $\lambda\le \frac{1}{2}$ and $s^{s+g+h}\ge s^{s+2g+1}\ge 2s^{s+2g}$, we have that  $2\lambda-1\le 0$ and $2s^{s+g+h}-4s^{s+2g}>0$, so the numerator of the second term is negative. Thus we can increase the value of this expression by choosing $B$ and $D$ as large as possible which yields:

$$\eff\left(A, B, C, D\right)\le\eff\left(A, \frac{\lambda}{1-\lambda}s^{s+g+h}, C, s^{s+g+h}\right).$$

We will do the same trick for $A$ and $C$.

$$\eff\left(A, \frac{\lambda}{1-\lambda}s^{s+g+h}, C, s^{s+g+h}\right)=$$
$$=2\frac{\lambda A+(1-\lambda)\frac{\lambda}{1-\lambda}s^{s+g+h} +\lambda C+(1-\lambda) s^{s+g+h} }{A+\frac{\lambda}{1-\lambda}s^{s+g+h}+C +s^{s+g+h}+2s^{s+g+h}-4s^{s+2g}} =$$
$$=2\frac{\lambda A+\lambda C+s^{s+g+h}}{A+C +(\frac{\lambda}{1-\lambda}+3)s^{s+g+h}-4s^{s+2g}} =$$
$$=2\lambda+2\frac{\left[1-\lambda(\frac{\lambda}{1-\lambda}+3)\right]s^{s+g+h}
+4\lambda s^{s+2g}}{A+C+\left(\frac{\lambda}{1-\lambda}+3\right)s^{s+g+h}-4s^{s+2g}}=$$
$$=2\lambda+2\frac{\frac{2\lambda^2-4\lambda+1}{1-\lambda}s^{s+g+h}
+4\lambda s^{s+2g}}{A+C+\left(\frac{\lambda}{1-\lambda}+3\right)s^{s+g+h}-4s^{s+2g}}.$$
If $\lambda=1-\frac{1}{\sqrt2}$, then $2\lambda^2-4\lambda+1=0$, so
$$\eff\left(A, \frac{\lambda}{1-\lambda}s^{s+g+h}, C, s^{s+g+h}\right)=$$
$$=2\lambda+2\frac{4\lambda s^{s+2g}}{A+C+\left(\frac{\lambda}{1-\lambda}+3\right)s^{s+g+h}-4s^{s+2g}}.$$

The numerator of the second term is $\Theta(s^{s+2g})$ while the denominator is $\Theta(s^{s+g+h})$ because $A\le \frac{1-\lambda}{\lambda}s^{s+g+h}$ and $C\le s^{s+g+h}$.
Thus
\[
\eff\left(A, \frac{\lambda}{1-\lambda}s^{s+g+h}, C, s^{s+g+h}\right)= 2\lambda+O(s^{g-h})=2-\sqrt{2}+o(1).
\qedm
\]

We have already proved that the number of those edges in the matching $M$ which are not in $M''$ is negligible. We can partition the rest of the edges (the regular ones) into disjoint subsets according to the $(a,b)$-bonded pairs determined by their beads. In every such $(a,b)$-bonded pair (for some $a<b$), we can repeat the argument of Lemma \ref{eff}. Hence, we can conclude that in this construction the size of a proper matching is at most $\left ( 2 - \sqrt 2+o(1) \right )n \approx 0.5858n$. \hfill $\Box$

\bigskip

\begin{rem}
One can similarly deal with the end of the proof of Lemma \ref{eff} in the heterogeneous case and then conclude that $\mu\le 2-\sqrt 2$, too.
\end{rem}

\section{Unbalanced necklaces} \label{unbalanced}

In this section, we consider the case when the number of red and blue beads are different. There are two possible measurements of unbalancedness we can use.  
On the one hand, we can define that the number of red beads is $\phi n$, while the number of blue beads is $n$  -- thus the total number of beads is $N=(1+\phi)n$.
On the other hand, we can first fix $N$, the total number of beads, and we define the number of red beads as $pN$, where $0< p< 1$.
(As before, we will omit integer parts).
Of course any of these parameters defines the other one:
\[
 p=\frac{\phi}{1+\phi},
\]
\[
\phi=\frac{p}{1-p}.
\]
In the calculations we will use $\phi$, however, we will show our final results both as a function of $\phi$ and as a function of $p$. We consider only the case $1/2\le \phi\le 2$, i.e., $1/3\le p\le 2/3$.

Now let $\mathcal{N}_\phi(N)$ denote the set of possible necklaces with $\phi n+n$ beads where $\phi n$ ones are red and $n$ ones are blue. Let $\mathcal{M}(L)$ denote the set of proper matchings for a given necklace $L$ in the heterogeneous model.
Moreover, let $\mu_\phi(N)=\min_{L\in\mathcal{N}_\phi(N)} \max_{M\in\mathcal{M}(L)} 2|M|$.  
We are interested in $\frac{\mu_\phi(N)}{N}$, the ratio of the covered beads to the total number of beads. Finally let $\mu_\phi=\limsup\limits_{N\rightarrow\infty} \mu_\phi(N)$. For the homogeneous model we similarly define $\mu^{\mathrm{hom}}_\phi(N)$ and $\mu^{\mathrm{hom}}_\phi$.

In this section we give a modified version of Construction \ref{constr_dirty}.
Let $t,u,v,w$ be real numbers such that $t+w=\phi,\; u+v=1,\; t\ge u$ and $v\ge w$. For the sake of simplicity, we will omit the floor and ceiling functions again even in the description of the construction.

\begin{constr}\label{constr_unbalanced}
Let $s\ge 2$ be a integer parameter, and let $n=s^{5s+1}$. The necklace consists of $s$ large arcs, each having $\phi s^{5s}$ red and $s^{5s}$ blue beads. Let $L_1,\dots,L_s$ denote the large arcs.

$L_i$ is divided into $s^{2s-i}$ red and $s^{2s-i}$ blue arcs, the colors alternates. Let $\ell_{i,j}$ denote the $j$\textsuperscript{th} arc of $L_i$, where $1\le j\le 2s^{2s-i}$. If $j$ is odd, then the arc $\ell_{i,j}$ always consists of $(t+u) s^{3s+i}$ beads.
If $j$ is even, then the arc $\ell_{i,j}$ always consists of $(v+w) s^{3s+i}$ beads.

In the next step, we will change the color of some beads in each $\ell_{i,j}$ in the following way.
If $j$ is odd, then let's divide each $\ell_{i,j}$ into $s^{s+2i}$ intervals of size $(t+u)s^{2s-i}$, and in each of this tiny intervals, change the color of $u s^{2s-i}$ beads backwards from the clockwise end of the tiny interval from red to blue.
If $j$ is even, then let's divide each $\ell_{i,j}$ into $s^{s+2i}$ intervals
of size $(v+w)s^{2s-i}$, and in each of this tiny intervals, change the color
of $w s^{2s-i}$ beads backwards from the clockwise end of the tiny interval
from blue to red. We will refer to those beads whose color were changed as dust in $\ell_{i,j}$.
\end{constr}

Notice that we cannot use the same formula for $a_k$, $b_k$, $c_k$, $d_k$ since there is no $\lambda$ in the current setup, but the definition of these parameters will be the same. Namely, these parameters will again denote the necessary number of ''small'' intervals (of size $O(s^{2s-h})$) to cover all the beads in the matching of the corresponding ''big'' intervals (of size $O(s^{2s-g})$) from the opposite side. Intentionally, we do not write down these definitions precisely because there would be too many very similar definitions depending on the color of the corresponding intervals. Otherwise, we will use the same notations. 
All estimations before the introduction of the efficiency notion 
in the analysis of Construction \ref{constr_dirty} work exactly the same way. 

In order to properly calculate the efficiency, we have to distinguish the homogeneous and heterogeneous case. 

\subsection{Homogeneous case}

The calculation will be different, if $\ell_{h,j}$ is blue, and if it is red. First, consider the blue case. We will use the same argument as in Lemma \ref{eff}.

Therefore, we can get that $\eff_{\blue} (\ell_{g,2i-1}\cup\ell_{g,2i})$ is at most 
$$\eff_{\blue}\left(A, B, C, D\right) =2\frac{wA+vB+wC+vD}{(v+w)(A+B+C+D)+(\phi+1)s^{s+g+h}+\Theta(s^{s+2g})},$$
where $A\le \frac{t}{w}s^{s+g+h}$, $B\le \frac{u}{v}s^{s+g+h}$, $C\le s^{s+g+h}$, $D\le s^{s+g+h}$.

Again, this expression is monotone increasing in $B$ and $D$, since $v\ge w$, so 
$$\eff_{\blue}\left(A, B, C, D\right)\le \frac{2}{v+w}\frac{w(A+C)+v(\frac{u}{v}+1)s^{s+g+h}}{A+C+(\frac{u}{v}+1)s^{s+g+h}+\frac{\phi+1}{v+w}s^{s+g+h}+\Theta(s^{s+2g})}=$$
$$=\frac{2w}{v+w}+\frac{2}{v+w}\frac{((v-w)(\frac{u}{v}+1)-\frac{\phi+1}{v+w}w)s^{s+g+h}+\Theta(s^{s+2g})}{A+C+(\frac{u}{v}+1)s^{s+g+h}+\frac{\phi+1}{v+w}s^{s+g+h}+\Theta(s^{s+2g})}.$$

Thus if we choose the parameters such that $(v-w)(\frac{u}{v}+1)-\frac{\phi+1}{v+w}w=0$, then we can conclude that $\eff_{\blue} (\ell_{g,2i-1}\cup\ell_{g,2i})\le \frac{2w}{v+w}+o(1)$.

Using the fact $u+v=1$, we get the quadratic equation $\left(\frac{v}{w}\right)^2-(\phi+1)\left(\frac{v}{w}\right)-1=0$ with one positive root
$$\eta_1=\frac{\phi+1+\sqrt{\left(\phi+1\right)^2+4}}{2},$$ thus $v=\eta_1 w$, and
$\eff_{\blue} (\ell_{g,2i-1}\cup\ell_{g,2i})\le \frac{2}{\eta_1+1}+o(1)$.

\bigskip

In the second case where $\ell_{h,j}$ is red, one can go through the similar arguments as above and can get to the point in which we choose the parameters such that the main term of the corresponding numerator is 0. Then by using the fact, $t+w=\phi$, we get that  $\left(\frac{t}{u}\right)^2-\frac{\phi+1}{\phi}\left(\frac{t}{u}\right)-1=0$, which has one positive root
$$\eta_2=\frac{\frac{\phi+1}{\phi}+\sqrt{\left(\frac{\phi+1}{\phi}\right)^2+4}}{2},$$ thus $t=\eta_2 u$, and
$\eff_{\red} (\ell_{g,2i-1}\cup\ell_{g,2i})\le \frac{2}{\eta_2+1}+o(1)$.

\bigskip

It is easy to check that the following system of equations has a unique solution. Also, if $\frac12\le \phi\le 2$, then $u, v, w, t\ge 0$.
$$t+w=\phi$$ $$u+v=1$$ $$v=\eta_1 w$$ $$t=\eta_2 u$$

Altogether we get that at most $\max\left(\frac{2}{\eta_1+1},\frac{2}{\eta_2+1}\right)+o(1)$ fraction of the beads are matched. By symmetry, we can assume that $\phi\ge 1$. Then $\eta_1\ge \eta_2$, thus $\frac{2}{\eta_1+1}\le \frac{2}{\eta_2+1}$, and
$$\frac{2}{\eta_2+1}=\frac{4}{\phi+3+\sqrt{\phi^2+2\phi+5}}.$$

One can get $\frac{\phi}{1+\phi}$ as a trivial lower bound (dashed green) for this case if we choose a secant which cut the red beads into two equal parts and consider a matching which covers all the red beads (and do not care about blue beads at all). In particular, for $\phi=2$ we have equality. (See Figures \ref{homogen} and \ref{homogen_p}.)

\begin{figure}[!ht]
\centering
\begin{tikzpicture}[line cap=round,line join=round,>=triangle 45,x=10.0cm,y=7.0cm]
\begin{axis}[
x=10.0cm,y=7.0cm,
axis lines=middle,
xmin=1.0,
xmax=2.0,
ymin=0.0,
ymax=1.0,
xtick={1.0,1.1,...,2.0},
ytick={0.0,0.1,...,1.2},
xlabel={$\varphi$}]
\clip(1.,0.) rectangle (2.,1.);
\draw[line width=2.pt,color=cblue,smooth,samples=100,domain=1.0:2.0] plot(\x,{4.0*(\x)/(3.0*(\x)+1.0+sqrt(5.0*(\x)^(2.0)+2.0*(\x)+1.0))});
\draw[line width=0.5pt,color=cred,smooth,samples=100,domain=1.0:2.0] plot(\x,{4.0/((\x)+3.0+sqrt(5.0+2.0*(\x)+(\x)^(2.0)))});

\draw[line width=0.5pt,color=cgreen, style=dashed, ,smooth,samples=100,domain=1.0:2.0] plot(\x,{(\x)/(1.0+(\x))});

\begin{scriptsize}
\draw[color=cblue] (1.9,0.75) node {$\frac{2}{\psi_2+1}$};
\draw[color=cred] (1.9,0.425) node {$\frac{2}{\psi_1+1}$};
\draw[color=cgreen] (1.9,0.6) node {$\frac{\varphi}{1+\varphi}$};

\draw (1.9,0.95) node[anchor=east] {\tiny Our best upper bound};
\draw[line width=2.pt,color=cblue] (1.9, 0.95) -- (2, 0.95);
\draw (1.9,0.9) node[anchor=east] {\tiny A trivial lower bound};
\draw[line width=0.5pt,color=cgreen, style=dashed] (1.9, 0.9) -- (2, 0.9);

\end{scriptsize}
\end{axis}
\end{tikzpicture}
\caption{Lower and upper bounds in the homogeneous case as a function of $\phi$, in the case $1\le \phi\le 2$.}
\label{homogen}
\end{figure}

\begin{figure}[!ht]
\centering
\begin{tikzpicture}[line cap=round,line join=round,>=triangle 45,x=30.0cm,y=7.0cm]
\begin{axis}[
x=30.0cm,y=7.0cm,
axis lines=middle,
xmin=0.333333333333,
xmax=0.666666666666,
ymin=0.0,
ymax=1.0,
xtick={0.30,0.35,...,0.700002},
ytick={-0.1,0.0,...,1.2},
xlabel={$p$}]]
%\clip(0.33333333,-0.16331276932332664) rectangle (0.666666666,1.2097268402855659);
\draw[line width=2.pt,color=cred,smooth,samples=100,domain=0.33333333:0.5] plot(\x,{4.0*(1.0-(\x))/((\x)+3.0*(1.0-(\x))+sqrt((\x)^(2.0)+2.0*(\x)*(1.0-(\x))+5.0*(1.0-(\x))^(2.0)))});
\draw[line width=0.5pt,color=cred,smooth,samples=100,domain=0.5:0.666666666] plot(\x,{4.0*(1.0-(\x))/((\x)+3.0*(1.0-(\x))+sqrt((\x)^(2.0)+2.0*(\x)*(1.0-(\x))+5.0*(1.0-(\x))^(2.0)))});
\draw[line width=0.5pt,color=cblue,smooth,samples=100,domain=0.33333333:0.5] plot(\x,{4.0*(\x)/(1.0-(\x)+3.0*(\x)+sqrt((1.0-(\x))^(2.0)+2.0*(\x)*(1.0-(\x))+5.0*(\x)^(2.0)))});
\draw[line width=2.pt,color=cblue,smooth,samples=100,domain=0.5:0.666666666] plot(\x,{4.0*(\x)/(1.0-(\x)+3.0*(\x)+sqrt((1.0-(\x))^(2.0)+2.0*(\x)*(1.0-(\x))+5.0*(\x)^(2.0)))});
\draw[line width=0.5pt,color=cgreen, style=dashed,smooth,samples=100,domain=0.33333333:0.5] plot(\x,{1-\x});
\draw[line width=0.5pt,color=cgreen, style=dashed, smooth,samples=100,domain=0.5:0.666666666] plot(\x,{\x});

\begin{scriptsize}
\draw[color=cblue] (1.9,0.75) node {$\frac{2}{\psi_2+1}$};
\draw[color=cred] (1.9,0.425) node {$\frac{2}{\psi_1+1}$};
\draw[color=cgreen] (1.9,0.6) node {$\frac{\varphi}{1+\varphi}$};

\draw (0.63,0.95) node[anchor=east] {\tiny Our best upper bound};
\draw[line width=2.pt,color=cred] (0.63, 0.95) -- (0.648, 0.95);
\draw[line width=2.pt,color=cblue] (0.648, 0.95) -- (0.666, 0.95);
\draw (0.63,0.9) node[anchor=east] {\tiny A trivial lower bound};
\draw[line width=0.5pt,color=cgreen, style=dashed] (0.63, 0.9) -- (0.666, 0.9);

\end{scriptsize}

\end{axis}
\end{tikzpicture}
\caption{Lower and upper bounds in the homogeneous case as a function of $p$, in the case $1/3\le p\le 2/3$.}
\label{homogen_p}
\end{figure}
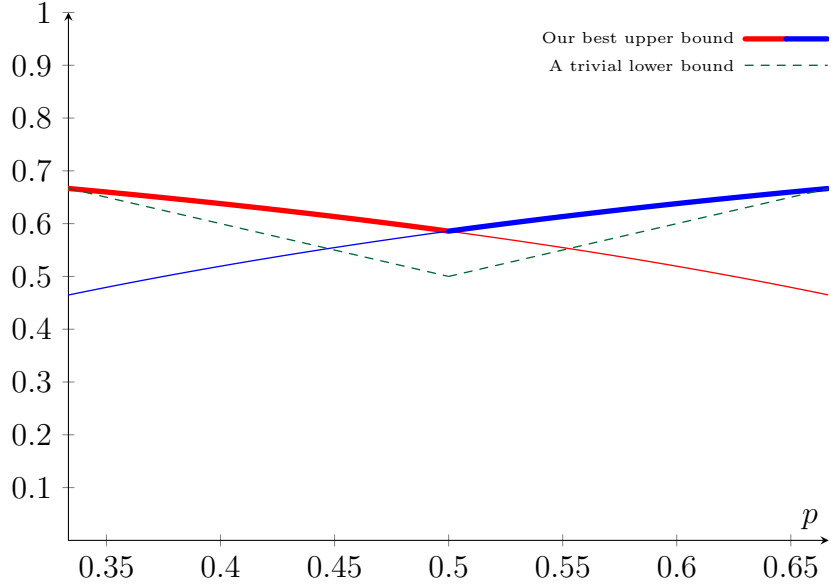

\subsection{Heterogeneous case}

The calculation will be different again depending on the color of $\ell_{h,j}$. First, consider the case when $\ell_{h,j}$ is blue. We will use the same argument as in Lemma \ref{eff}.

Therefore, we can get that $\eff_{\blue} (\ell_{g,2i-1}\cup\ell_{g,2i})$ is at most 
$$\eff_{\blue}\left(A, B, C, D\right) =2\frac{vA+wB+vC+wD}{(v+w)(A+B+C+D)+(\phi+1)s^{s+g+h}+\Theta(s^{s+2g})},$$
where $A\le \frac{t}{v}s^{s+g+h}$, $B\le \frac{u}{w}s^{s+g+h}$, $C\le \frac{w}{v}s^{s+g+h}$, $D\le \frac{v}{w}s^{s+g+h}$.

Again, this expression is monotone increasing in $A$ and $C$, since $v\ge w$, so 
$$\eff_{\blue}\left(A, B, C, D\right)\le \frac{2}{v+w}\frac{w(B+D)+v(\frac{t}{v}+\frac{w}{v})s^{s+g+h}}{B+D+(\frac{t}{v}+\frac{w}{v})s^{s+g+h}+\frac{\phi+1}{v+w}s^{s+g+h}+\Theta(s^{s+2g})}=$$
$$=\frac{2w}{v+w}+\frac{2}{v+w}\frac{((v-w)(\frac{t}{v}+\frac{w}{v})-\frac{\phi+1}{v+w}w)s^{s+g+h}+\Theta(s^{s+2g})}{B+D+(\frac{u}{v}+1)s^{s+g+h}+\frac{\phi+1}{v+w}s^{s+g+h}+\Theta(s^{s+2g})}.$$

Thus if we choose the parameters such that $(v-w)(\frac{t}{v}+\frac{w}{v})-\frac{\phi+1}{v+w}w=0$, then we can conclude that $\eff_{\blue} (\ell_{g,2i-1}\cup\ell_{g,2i})\le \frac{2w}{v+w}+o(1)$.

Using the fact, $t+w=\phi$, we get the quadratic equation $\left(\frac{v}{w}\right)^2-\frac{\phi+1}{\phi}\left(\frac{v}{w}\right)-1=0$ with one positive root which we called $\eta_2$, and
$\eff_{\blue} (\ell_{g,2i-1}\cup\ell_{g,2i})\le \frac{2}{\eta_2+1}+o(1)$.

If $\ell_{h,j}$ is red, then we get that
$\left(\frac{t}{u}\right)^2-(\phi+1)      \left(\frac{t}{u}\right)-1=0$, which has one positive root $\eta_1$, and
$\eff_{\red} (\ell_{g,2i-1}\cup\ell_{g,2i})\le \frac{2}{\eta_1+1}+o(1)$.

\bigskip

It is easy to check that the following system of equations has a unique solution, and $u, v, w, t$ are always positive.
$$t+w=\phi$$ $$u+v=1$$ $$v=\eta_2 w$$ $$t=\eta_1 u$$

Altogether we get that at most $\max\left(\frac{2}{\eta_1+1},\frac{2}{\eta_2+1}\right)+o(1)$ fraction of the beads are matched, which is the same as in the homogeneous case. By symmetry, we can assume that $\phi\ge 1$. Then $\eta_1\ge \eta_2$, thus $\frac{2}{\eta_1+1}\le \frac{2}{\eta_2+1}$, and
$$\frac{2}{\eta_2+1}=\frac{4}{\phi+3+\sqrt{\phi^2+2\phi+5}}.$$

One can get $\frac{1}{1+\phi}$ as a trivial lower bound (dashed green) for
this case if we choose a secant which cut the blue beads into two equal
parts. Since $\phi\ge 1$, one can find a matching which covers all the blue
beads from the side which contains at most as many red beads as the other
side. On the other hand, one can get $\frac{2}{1+\phi}$ as a trivial upper
bound (dotted yellow) because in every edge of the matching exactly one of the
endpoints have color blue.  In particular, for $\phi=2$ we have equality with
the trivial bound, thus our construction gives a non-trivial bound only for
$\phi<2$. (See Figures  \ref{heterogen} and \ref{heterogen_p}.)

\begin{figure}[!ht]
\centering
\definecolor{ffxfqq}{rgb}{1.,0.4980392156862745,0.}
\definecolor{ffqqqq}{rgb}{1.,0.,0.}
\definecolor{qqwuqq}{rgb}{0.,0.39215686274509803,0.}
\definecolor{qqqqff}{rgb}{0.,0.,1.}
\begin{tikzpicture}[line cap=round,line join=round,>=triangle 45,x=10.0cm,y=7.0cm]
\begin{axis}[
x=10.0cm,y=7.0cm,
axis lines=middle,
xmin=1.0,
xmax=2.0,
ymin=0.0,
ymax=1.0,
xtick={1.0,1.1,...,2.0},
ytick={0.0,0.1,...,1.2},
xlabel={$\varphi$}]
\clip(1.,0.) rectangle (2.,1.);
\draw[line width=2.pt,color=cblue,smooth,samples=100,domain=1.0:2.0] plot(\x,{4.0*(\x)/(3.0*(\x)+1.0+sqrt(5.0*(\x)^(2.0)+2.0*(\x)+1.0))});
\draw[line width=0.5pt,color=cred,smooth,samples=100,domain=1.0:2.0] plot(\x,{4.0/((\x)+3.0+sqrt(5.0+2.0*(\x)+(\x)^(2.0)))});

\draw[line width=1pt,color=cyellow,style=dotted, smooth,samples=100,domain=1.0:2.0] plot(\x,{(2)/(1.0+(\x))});
\draw[line width=1pt,color=cgreen, style=dashed, smooth,samples=100,domain=1.0:2.0] plot(\x,{1.0/(1.0+(\x))});

\begin{scriptsize}
\draw[color=cblue] (1.1,0.65) node {$\frac{2}{\psi_2+1}$};
\draw[color=cred] (1.1,0.525) node {$\frac{2}{\psi_1+1}$};
\draw[color=cyellow] (1.1,0.9) node {$\frac{2}{1+\varphi}$};
\draw[color=cgreen] (1.1,0.4) node {$\frac{1}{1+\varphi}$};

\draw (1.9,0.95) node[anchor=east] {\tiny Our best upper bound};
\draw[line width=2.pt,color=cblue] (1.9, 0.95) -- (2, 0.95);
\draw (1.9,0.9) node[anchor=east] {\tiny A trivial upper bound};
\draw[line width=1pt,color=cyellow, style=dotted] (1.9, 0.9) -- (2, 0.9);
\draw (1.9,0.85) node[anchor=east] {\tiny A trivial lower bound};
\draw[line width=0.5pt,color=cgreen, style=dashed] (1.9, 0.85) -- (2, 0.85);

\end{scriptsize}
\end{axis}
\end{tikzpicture}
\caption{Lower and upper bounds in the heterogeneous case
 as a function of $\phi$, in the case $1\le \phi\le 2$.}
\label{heterogen}
\end{figure}

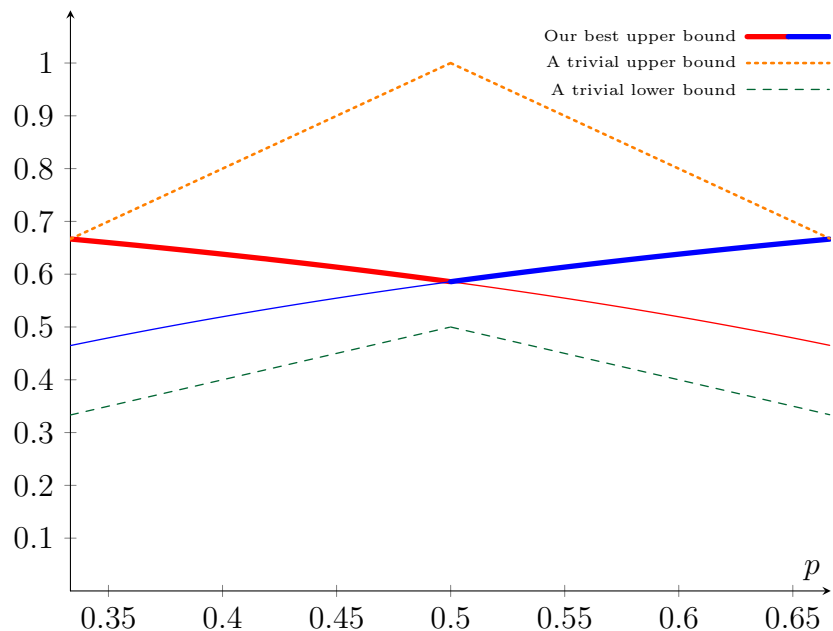
\begin{figure}[!ht]
\centering
\definecolor{ffxfqq}{rgb}{1.,0.4980392156862745,0.}
\definecolor{ffqqqq}{rgb}{1.,0.,0.}
\definecolor{qqwuqq}{rgb}{0.,0.39215686274509803,0.}
\definecolor{qqqqff}{rgb}{0.,0.,1.}
\begin{tikzpicture}[line cap=round,line join=round,>=triangle 45,x=30.0cm,y=7.0cm]
\begin{axis}[
x=30.0cm,y=7.0cm,
axis lines=middle,
xmin=0.333333333333,
xmax=0.666666666666,
ymin=0.0,
ymax=1.1,
xtick={0.30,0.35,...,0.700002},
ytick={-0.1,0.0,...,1.1},
xlabel={$p$}]]
%\clip(0.33333333,-0.16331276932332664) rectangle (0.666666666,1.2097268402855659);
\draw[line width=2.pt,color=cred,smooth,samples=100,domain=0.33333333:0.5] plot(\x,{4.0*(1.0-(\x))/((\x)+3.0*(1.0-(\x))+sqrt((\x)^(2.0)+2.0*(\x)*(1.0-(\x))+5.0*(1.0-(\x))^(2.0)))});
\draw[line width=0.5pt,color=cred,smooth,samples=100,domain=0.5:0.666666666] plot(\x,{4.0*(1.0-(\x))/((\x)+3.0*(1.0-(\x))+sqrt((\x)^(2.0)+2.0*(\x)*(1.0-(\x))+5.0*(1.0-(\x))^(2.0)))});
\draw[line width=0.5pt,color=cblue,smooth,samples=100,domain=0.33333333:0.5] plot(\x,{4.0*(\x)/(1.0-(\x)+3.0*(\x)+sqrt((1.0-(\x))^(2.0)+2.0*(\x)*(1.0-(\x))+5.0*(\x)^(2.0)))});
\draw[line width=2.pt,color=cblue,smooth,samples=100,domain=0.5:0.666666666] plot(\x,{4.0*(\x)/(1.0-(\x)+3.0*(\x)+sqrt((1.0-(\x))^(2.0)+2.0*(\x)*(1.0-(\x))+5.0*(\x)^(2.0)))});
\draw[line width=1pt,color=cyellow, style=dotted,smooth,samples=100,domain=0.33333333:0.5] plot(\x,{2*\x});
\draw[line width=1pt,color=cyellow,style=dotted, smooth,samples=100,domain=0.5:0.666666666] plot(\x,{2-2*\x});
\draw[line width=0.5pt,color=cgreen,style=dashed, smooth,samples=100,domain=0.33333333:0.5] plot(\x,{\x});
\draw[line width=0.5pt,color=cgreen,style=dashed, smooth,samples=100,domain=0.5:0.666666666] plot(\x,{1-\x});

\draw[color=cblue] (1.9,0.75) node {$\frac{2}{\psi_2+1}$};
\draw[color=cred] (1.9,0.425) node {$\frac{2}{\psi_1+1}$};
\draw[color=cgreen] (1.9,0.6) node {$\frac{\varphi}{1+\varphi}$};

\draw (0.63,1.05) node[anchor=east] {\tiny Our best upper bound};
\draw[line width=2.pt,color=cred] (0.63, 1.05) -- (0.648, 1.05);
\draw[line width=2.pt,color=cblue] (0.648, 1.05) -- (0.666, 1.05);
\draw (0.63,1.0) node[anchor=east] {\tiny A trivial upper bound};
\draw[line width=1 pt,color=cyellow, style=dotted] (0.63, 1.0) -- (0.666, 1.0);
\draw (0.63,0.95) node[anchor=east] {\tiny A trivial lower bound};
\draw[line width=0.5pt,color=cgreen, style=dashed] (0.63, 0.95) -- (0.666, 0.95);

\end{axis}
\end{tikzpicture}
\caption{Lower and upper bounds in the heterogeneous case
 as a function of $p$, in the case $1/3\le p\le 2/3$.}
\label{heterogen_p}
\end{figure}

\section*{Acknowledgement}
This research was started during the 9th Emléktábla Workshop, 2019.
The authors are thankful for the organizers for inviting them.
We also thank D.~P\'alv\"olgyi, G.~Damásdi, T.~Fleiner and Zs.~Jank\'o
for valuable questions and observations.


\begin{thebibliography}{1}

\bibitem{LP} R.~B.~Lyngs{\o} and C.~N.~S.~Pedersen,
Protein Folding in the 2D HP Model,
\href{https://tidsskrift.dk/brics/article/view/20073}
%, doi: 10.7146/brics.v6i16.20073.
{\newblock {\em  BRICS Report Series}, {\bf RS-99-16}, (1999).}
 
\bibitem{KPT1} J. Kyn\v cl, J. Pach and G. T\'oth, Long alternating paths in bicolored point sets, 
\newblock in Graph Drawing (J. Pach, ed.), Lecture Notes in Computer
Science {\bf 3383}, Springer-Verlag, Berlin, (2004), pp.\ 340--348. 

\bibitem{KPT2} J. Kyn\v cl, J. Pach and G. T\'oth, Long alternating paths in bicolored point sets, 
\newblock {\em Discrete Mathematics}, {\bf 308}, (2008), pp.\ 4315--4322.

\bibitem{Abellanas} M.~Abellanas, A.~Garcia, F.~Hurtado and J.~Tejel, Caminos alternantes, in: {\em X Encuentros de Geometria Computacional (in Spanish)}, Sevilla, (2003), pp.\ 7--12.

\bibitem{Sebo} G.~Brevier, M.~Preissmann and A.~Seb\H{o},
personal communication (2004). 

\bibitem{HM} P.~Hajnal and V.~M\'esz\'aros, A note on noncrossing path in colored convex sets, \newblock {manuscript, (2010)}.

\bibitem{Viola} V.~M\'esz\'aros, Extremal problems on planar point sets,
\newblock {\em Ph.D. thesis}, \href{http://doktori.bibl.u-szeged.hu/688/1/mvdoktori.pdf}{doktori.bibl.u-szeged.hu/688/1/mvdoktori.pdf}, (2011).

\bibitem{Viola2} V.~M\'esz\'aros, Separated matchings and small discrepancy colorings. \newblock{\em Computational Geometry}, {\em Lecture Notes in Comput. Sci.}, {\bf 7579}, Springer, Cham, (2011), pp.\ 236--248.

\bibitem{MR} C.~M\"ullner and A.~Ryzhikov,
Palindromic Subsequences in Finite Words.
\newblock \href{https://arxiv.org/abs/1901.07502}
{{\em  arXiv}, {\bf 1901.07502}, (2019).}

\bibitem{MR2} C.~M\"ullner and A.~Ryzhikov, 
Palindromic subsequences in finite words. 
\newblock In Proc
13th Int.\ Conf. Language and Automata Theory and Applications (LATA), (2019), pp.\ 460-–468.% doi:10.1007/978-3-030-13435-8\_34.

\bibitem{MV} W. Mulzer and P. Valtr, Long alternating paths exist. \newblock \href{https://arxiv.org/abs/2003.13291}
{{\em arXiv}, {\bf 2003.13291}, (2020).}


\end{thebibliography}
\end{document}